\documentclass{article}
\usepackage{amssymb}
\usepackage{amsmath}
\usepackage{amsfonts}
\usepackage{amscd}
\usepackage[dvips]{graphicx}

\newtheorem{defi}{D\'{e}finition}
\newtheorem{propo}{Proposition}

\newtheorem{theorem}{Theorem}
\newtheorem{definition}{Definition}

\begin{document}

\title{Linear Algebra in the vector space of intervals $\overline{\mathbb{IR}%
}$}
\author{Nicolas GOZE }
\date{}
\maketitle

\begin{abstract}
In a previous paper, we have given an algebraic model to the set of intervals. Here, we apply this model in a linear frame.
We define a notion of diagonalization of square matrices whose coefficients are intervals. But in this case, with respect to the real case,
a matrix of order $n$ could have more than $n$ eigenvalues (the set of intervals is not factorial). We consider a notion of
central eigenvalues permits to describe criterium of diagonalization. As application, we define a notion of Exponential mapping.
\end{abstract}
\section{The associative algebra $\overline{\mathbb{IR}}$}

In \cite{GN-R}, we have given a representation of the set of intervals in
terms of associative algebra. More precisely, we define on the set $\mathbb{%
IR}$ of intervals of $\mathbb{R}$ a $\mathbb{R}$-vector space structure.
Next we embed $\mathbb{IR}$ in a $4$-dimensional associative algebra. This
embedding permits to describe a unique distributive multiplication which
contains all the possible results of the usual product of intervals and the
monotony property is always conserved. Moreover, this new product is minimal with
respect the distributivity and the monotony properties.

\noindent In this section, we present briefly this construction (for more
details, see \cite{GN-R}).

\subsection{Vector space structure on $\mathbb{IR}$}

Let $\mathbb{IR}$ be the set of intervals of $\mathbb{R}$, that is
\begin{equation*}
\mathbb{IR=\{[}a,b],a,b\in \mathbb{R}\text{, }a\leq b\}.
\end{equation*}%
This set is provided with a semi-group structure that we can complete as
follow: we consider the equivalence relation on $\mathbb{IR}\times \mathbb{IR%
}$:
\begin{equation*}
(X,Y)\sim (Z,T)\Longleftrightarrow X+T=Y+Z
\end{equation*}%
for all $X,Y,Z,T\in \mathbb{IR}$. The quotient set is denoted by $\overline{%
\mathbb{IR}}$. The addition of intervals is compatible with this equivalence
relation :
\begin{equation*}
\overline{(X,Y)}+\overline{(z,t)}=\overline{(x+z,y+t)}
\end{equation*}%
where $\overline{(x,y)}$ is the equivalence class of $(X,Y).$ The unit is $%
\overline{0}=\{(X,X),X\in \mathbb{IR\}}$ and each element has an inverse
\begin{equation*}
\smallsetminus \overline{(X,Y)}=\overline{(Y,X)}.
\end{equation*}%
Then $(\overline{\mathbb{IR}},+)$ is a commutative group. We prove also in
\cite{GN-R}, that any equivalence class admits a canonical representant of
type $(K,0)$ or $(0,K)$ or $(a,a)$, with $K\in \mathbb{IR}$ and $a=[a,a],\
a\in \mathbb{R}$. We provides the group $(\overline{\mathbb{IR}},+)$ with a
real vector space structure, the external product being given by
\begin{equation*}
\left\{
\begin{array}{c}
\alpha \cdot \overline{(K,0)}\text{\ }=\overline{(\alpha K,0)}\text{ } \\
\alpha \cdot \overline{(0,K)}\text{\ }=\overline{(0,\alpha K)}\text{ }%
\end{array}%
\right.
\end{equation*}%
if $\alpha >0.$ If $\alpha <0$ we put $\beta =-\alpha $ and
\begin{equation*}
\left\{
\begin{array}{c}
\alpha \cdot \overline{(K,0)}\text{\ }=\overline{(0,\beta K)}\text{ } \\
\alpha \cdot \overline{(0,K)}\text{\ }=\overline{(\beta K,0)}\text{ }%
\end{array}%
\right. .
\end{equation*}%
The triplet $(\overline{\mathbb{IR}},+,\cdot )$ is a real vector space.

\medskip

\noindent\textbf{Remark.} To simplify notations, we write $(K,0)$ or $(0,K)$
in place of $\overline{(K,0)}$ or $\overline{(0,K)}$.

\subsection{The associative algebra $\mathcal{A}_4$}

Recall that by an algebra we mean a real vector space with an associative
ring structure. Consider the $4$-dimensional associative algebra whose
product in a basis $\{e_{1},e_{2},e_{3},e_{4}\}$ is given by
\begin{equation*}
\begin{tabular}{|l|l|l|l|l|}
\hline
& $e_{1}$ & $e_{2}$ & $e_{3}$ & $e_{4}$ \\ \hline
$e_{1}$ & $e_{1}$ & $0$ & $0$ & $e_{4}$ \\ \hline
$e_{2}$ & $0$ & $e_{2}$ & $e_{3}$ & $0$ \\ \hline
$e_{3}$ & $0$ & $e_{3}$ & $e_{2}$ & $0$ \\ \hline
$e_{4}$ & $e_{4}$ & $0$ & $0$ & $e_{1}$ \\ \hline
\end{tabular}%
.
\end{equation*}%
The unit is the vector $e_{1}+e_{2}.$ This algebra is a direct sum of two
ideals: $\mathcal{A}_{4}=I_{1}+I_{2}$ where $I_{1}$ is generated by $e_{1}$
and $e_{4}$ and $I_{2}$ is generated by $e_{2}$ and $e_{3}.$ It is not an
integral domain, that is, we have divisors of $0.$ For example $e_{1}\cdot
e_{2}=0.$ The ring $\mathcal{A}_{4}$ is principal that is every ideal is
generated by one element. The cartesian expression of this product is, for $%
x=(x_{1},x_{2},x_{3},x_{4})$ and $y=(y_{1},y_{2},y_{3},y_{4})$ in $\mathcal{A%
}_{4}$:
\begin{equation*}
x\cdot
y=(x_{1}y_{1}+x_{4}y_{4},x_{2}y_{2}+x_{3}y_{3},x_{3}y_{2}+x_{2}y_{3},x_{4}y_{1}+x_{1}y_{4}).
\end{equation*}%
The multiplicative group $\mathcal{A}_{4}^{\ast }$ \ of invertible elements
is the set of elements $x=(x_{1},x_{2},x_{3},x_{4})$ such that
\begin{equation*}
\left\{
\begin{array}{c}
x_{4}\neq \pm x_{1}, \\
x_{3}\neq \pm x_{2}.%
\end{array}%
\right.
\end{equation*}%
If $x\in $ $\mathcal{A}_{4}^{\ast }$ we have:
\begin{equation*}
x^{-1}=\left( \frac{x_{1}}{x_{1}^{2}-x_{4}^{2}},\frac{x_{2}}{%
x_{2}^{2}-x_{3}^{2}},\frac{x_{3}}{x_{2}^{2}-x_{3}^{2}},\frac{x_{4}}{%
x_{1}^{2}-x_{4}^{2}}\right) .
\end{equation*}

\subsection{A product in $\overline{\mathbb{IR}}$}

We define a correspondence between $\mathbb{IR}$ and $\mathcal{A}$. Let $%
\varphi$ be the map
\begin{equation*}
\varphi : \overline{\mathbb{IR}} \longrightarrow \mathcal{A}_{4}
\end{equation*}%
given by
\begin{equation*}
\varphi (K,0)=\left\{
\begin{array}{l}
(x_{1},x_{2},0,0)\text{ if }x_{1},x_{2}\geq 0, \\
(0,x_{2},-x_{1},0)\text{ if }x_{1}\leq 0\text{ and }x_{2}\geq 0, \\
(0,0,-x_{1},-x_{2})\text{ if }x_{1},x_{2}\leq 0,%
\end{array}%
\right.
\end{equation*}
and
\begin{equation*}
\overline{\varphi }(0,K)=-\overline{\varphi }(K,0).
\end{equation*}%
for any $K=[x_{1},x_{2}]$. Then we have
\begin{equation*}
\begin{tabular}{lll}
$\overline{\varphi }(0,K)$ & $=$ & $\left\{
\begin{array}{l}
-x_{1},-x_{2},0,0)\text{ if }x_{1}\geq 0, \\
(0,-x_{2},-x_{1},0)\text{ if }x_{1}x_{2}\leq 0, \\
(0,0,-x_{1},-x_{2})\text{ if }x_{2}\leq 0.%
\end{array}%
\right. $%
\end{tabular}%
\end{equation*}%
Thus the image of $\overline{\mathbb{IR}}$ in $\mathcal{A}_{4}$ is
constituted of the elements
\begin{equation*}
\left\{
\begin{array}{l}
(x_{1},x_{2},0,0)\text{ with }0\leq x_{1}\leq x_{2}\text{ which corresponds
to }([x_{1},x_{2}],0), \\
(0,x_{2},-x_{1},0)\text{ with }x_{1}\leq 0\leq x_{2}\text{ which corresponds
to }([x_{1},x_{2}],0), \\
(0,0,-x_{1},-x_{2})\text{ with }x_{1}\leq x_{2}\leq 0\text{ which
corresponds to }([x_{1},x_{2}],0), \\
(-x_{1},-x_{2},0,0)\text{with }0\leq x_{1}\leq x_{2}\text{ which corresponds
to }(0,[x_{1},x_{2}]), \\
(0,-x_{2},x_{1},0)\text{ with }x_{1}\leq 0\leq x_{2}\text{ which corresponds
to }(0,[x_{1},x_{2}]), \\
(0,0,x_{1},x_{2})\text{ with }x_{1}\leq x_{2}\leq 0\text{ which corresponds
to }(0,[x_{1},x_{2}]).%
\end{array}%
\right.
\end{equation*}%
The map $\overline{\varphi }:$ $\overline{\mathbb{IR}}\mathbb{%
\longrightarrow }\mathcal{A}_{4}$ \ is not linear. We introduce in $\mathcal{%
A}_{4}$ the following equivalence relation $\mathcal{R}$ given by
\begin{equation*}
(x_{1},x_{2},x_{3},x_{4})\sim (y_{1},y_{2},y_{3},y_{4})\Longleftrightarrow
\left\{
\begin{array}{c}
x_{1}-y_{1}=x_{3}-y_{3}, \\
x_{2}-y_{2}=x_{4}-y_{4}%
\end{array}%
\right.
\end{equation*}%
and consider the map
\begin{equation*}
\overline{\overline{\varphi }}:\overline{\mathbb{IR}}\longrightarrow
\overline{\mathcal{A}_{4}}=\frac{\mathcal{A}_{4}}{\mathcal{R}}
\end{equation*}%
given by $\overline{\overline{\varphi }}=\Pi \circ \overline{\varphi }$
where $\Pi $ is a canonical projection. This map is surjective. In fact we
have the correspondence

\begin{itemize}
\item $x_{1}-x_{3}\geq 0$, $x_{2}-x_{4}\geq 0$, $x_{1}-x_{3}\leq x_{2}-x_{4}$%
\begin{equation*}
(x_{1},x_{2},x_{3},x_{4})\sim (x_{1}-x_{3},x_{2}-x_{4},0,0)=\overline{%
\varphi }([x_{1}-x_{3},x_{2}-x_{4}],0).
\end{equation*}

\item $x_{1}-x_{3}\geq 0$, $x_{2}-x_{4}\geq 0$, $x_{1}-x_{3}\geq x_{2}-x_{4}$%
\begin{equation*}
(x_{1},x_{2},x_{3},x_{4})\sim (0,0,x_{3}-x_{1},x_{4}-x_{2})=\overline{%
\varphi }(0,[x_{3}-x_{1},x_{4}-x_{2}]).
\end{equation*}

\item $x_{1}-x_{3}\geq 0$, $x_{2}-x_{4}\leq 0$%
\begin{equation*}
(x_{1},x_{2},x_{3},x_{4})\sim (0,x_{2}-x_{4},x_{3}-x_{1},,0)=\overline{%
\varphi }(0,[x_{3}-x_{1},x_{4}-x_{2}]).
\end{equation*}

\item $x_{1}-x_{3}\leq 0$, $x_{2}-x_{4}\geq 0$%
\begin{equation*}
(x_{1},x_{2},x_{3},x_{4})\sim (0,x_{2}-x_{4},x_{3}-x_{1},,0)=\overline{%
\varphi }([x_{3}-x_{1},x_{2}-x_{4}],0).
\end{equation*}

\item $x_{1}-x_{3}\leq 0$, $x_{2}-x_{4}\leq 0$, $x_{1}-x_{3}\geq x_{2}-x_{4}$%
\begin{equation*}
(x_{1},x_{2},x_{3},x_{4})\sim (x_{1}-x_{3},x_{2}-x_{4},0,0)=\overline{%
\varphi }(0,[x_{3}-x_{1},x_{4}-x_{2}]).
\end{equation*}

\item $x_{1}-x_{3}\leq 0$, $x_{2}-x_{4}\leq 0$, $x_{1}-x_{3}\leq
x_{2}-x_{4}(x_{1},x_{2},x_{3},x_{4})\sim (0,0,x_{3}-x_{1},x_{4}-x_{2})=%
\overline{\varphi }([x_{3}-x_{1},x_{2}-x_{4}],0).$
\end{itemize}

\noindent This correspondence defines a map%
\begin{equation*}
\psi :\overline{\mathcal{A}_{4}}\longrightarrow \overline{\mathbb{IR}}.
\end{equation*}%
In the following, to simplify notation, we write $\overline{\varphi }$
instead of $\overline{\overline{\varphi }}.$

\begin{definition}
For any $\mathcal{X},\mathcal{X}^{\prime }\in \overline{\mathbb{IR}}$, we
put
\begin{equation*}
\mathcal{X}\bullet \mathcal{X}^{\prime }=\psi (\overline{\overline{\varphi }(%
\mathcal{X})}\bullet \overline{\overline{\varphi }(\mathcal{X}^{\prime })}).
\end{equation*}
\end{definition}

This multiplication is distributive with respect the the addition. In fact
\begin{equation*}
(\mathcal{X}_{1}+\mathcal{X}_{2})\bullet \mathcal{X}^{\prime }=\psi (%
\overline{\overline{\varphi }(\mathcal{X}_{1}+\mathcal{X}_{2})}\bullet
\overline{\overline{\varphi }(\mathcal{X}^{\prime })}).
\end{equation*}%
Suppose that $\overline{\varphi }(\mathcal{X}_{1}+\mathcal{X}_{2})\neq
\overline{\varphi }(\mathcal{X}_{1})+\overline{\varphi }(\mathcal{X}_{2}).$
In this case this means that $\overline{\varphi }(\mathcal{X}_{1})+\overline{%
\varphi }(\mathcal{X}_{2})\notin \ Im \overline{\varphi }.$ But by
construction $\overline{\overline{\varphi }(\mathcal{X}_{1}+\mathcal{X}_{2})}%
\in \ Im \overline{\varphi }$ and this coincides with $\overline{\varphi }(%
\mathcal{X}_{1}+\mathcal{X}_{2}).$ For numerical application of this product, see (\cite[GN-R]).

\section{The module $gl(n,\overline{\mathbb{IR}})$}

Let $gl(n,\overline{\mathbb{IR}})$ be the set of square matrices of order $n$
whose elements are in $\overline{\mathbb{IR}}.$ A matrice of $gl(n,\overline{\mathbb{IR}})$ is  denoted by%
\begin{equation*}
A=\left( \mathcal{X}_{ij}\right) _{i,j=1,\cdots ,n}
\end{equation*}%
with $\mathcal{X}_{ij}=(K_{ij},0)$ or $(0,K_{ij}).$ It is clear that $gl(n,%
\overline{\mathbb{IR}})$ is a real vector space. We define a product on it
puting
\begin{equation*}
A\cdot B=\left( \mathcal{X}_{ij}\right) \cdot \left( Y_{ij}\right) =\left(
Z_{ij}\right)
\end{equation*}%
with $Z_{ij}=\sum\limits_{k=1}^{n}\mathcal{X}_{ik}\cdot Y_{ik}.$ This last
product being the associative product on $\overline{\mathbb{IR}}.$ Thus $%
gl(n,\overline{\mathbb{IR}})$ is an associative algebra.

\begin{defi}
A matrice $A\in gl(n,\overline{\mathbb{IR}})$ is called inversible if its
determinant, computed by the Cramer rule, is an inversible element in $\overline{\mathbb{IR}}$.
\end{defi}

\noindent Recall that the group $\overline{\mathbb{IR}}$ of
inversible elements contain
\begin{equation*}
\mathcal{X}_{i}=(K_{i},0)\text{ or }(0,K_{i})
\end{equation*}%
with $0\notin K_{i}.$ To compute the determinant, we use the classical
formula of Cramer.

\bigskip

\noindent \textbf{Example 1. }Let us consider the matrix%
\begin{equation*}
M=\left(
\begin{array}{cc}
\lbrack 1,2] & [-1.3] \\
\lbrack -1,3] & [1,2]%
\end{array}%
\right)
\end{equation*}%
Thus
\begin{eqnarray*}
\det B_1 &=&([1,2],0)([1,2],0)\smallsetminus ([-1,3],0)([-1,3],0) \\
&=&([1,4],0)\smallsetminus ([-3,9],0) \\
&=&([0,[-4,5]) \\
&=&\smallsetminus ([-4,5],0).
\end{eqnarray*}%
As $([-4,5],0)$ is not an inversible element of $\overline{\mathbb{IR}},$
the matrix $B_1$ is not inversible.

\bigskip

\noindent \textbf{Example 2. }Now if
\begin{equation*}
B_2=\left(
\begin{array}{cc}
\lbrack 1,2] & [-1.3] \\
\lbrack -1,3] & [1,7]%
\end{array}%
\right)
\end{equation*}%
then, by the similar computation, we obtain%
\begin{equation*}
\det B_2=\smallsetminus ([-7,-4],0)
\end{equation*}%
and $B_2$ is invertible.

\medskip

\begin{definition}
If $A$ is an invertible matrix on $gl(n,\overline{\mathbb{IR}})$, the
inverse matrix $A^{-1}$ of $A$ is given by
\begin{equation*}
A\cdot A^{-1}=Id
\end{equation*}
where
\begin{equation*}
Id=\left(
\begin{array}{cccc}
1 & 0 & \cdots & 0 \\
0 & 1 & \cdots & 0 \\
\vdots & \vdots & \ddots & \vdots \\
0 & 0 & \cdots & 1%
\end{array}%
\right)
\end{equation*}%
with $1=([1,1],0)$ and $0=([0,0],0).$
\end{definition}

The determination of $A^{-1}$ can be computed using the classical rules.

\medskip

\noindent \textbf{Example. }If we consider the invertible matrix $B_2$, we
obtain
\begin{equation*}
B_2^{-1}=[\frac{1}{7},\frac{1}{4}]\left(
\begin{array}{cc}
\lbrack 1,7] & \smallsetminus \lbrack -1.3] \\
\smallsetminus \lbrack -1,3] & [1,2]%
\end{array}%
\right) .
\end{equation*}%
Let us verify that $B_2B_2^{-1}=Id$. Using the product on $\overline{\mathbb{%
IR}} $ we obtain
\begin{equation*}
B_2B_2^{-1}=[\frac{1}{7},\frac{1}{4}]\left(
\begin{array}{cc}
\lbrack 1,2] & [-1.3] \\
\lbrack -1,3] & [1,7]%
\end{array}%
\right) \cdot \left(
\begin{array}{cc}
\lbrack 1,7] & \smallsetminus \lbrack -1.3] \\
\smallsetminus \lbrack -1,3] & [1,2]%
\end{array}%
\right) .
\end{equation*}%
The coefficient in place $(1,1)$ is
\begin{equation*}
a_{11}=[\frac{1}{7},\frac{1}{4}]([1,2][1,7]+[-1.3](\smallsetminus \lbrack
-1,3])).
\end{equation*}%
>From the definition of the product (see section 1), this element is
\begin{eqnarray*}
a_{11} &=&(\frac{1}{7},\frac{1}{4},0,0)((1,2,0,0)(1,7,0,0)-(0,3,1,0)(0,3,1,0)
\\
&=&(\frac{1}{7},\frac{1}{4},0,0)((1,14,0,0)-(0,10,6,0)) \\
&=&(\frac{1}{7},\frac{1}{4},0,0)(1,4,-6,0) \\
&=&(\frac{1}{7},\frac{1}{4},0,0)(7,4,0,0) \\
&=&(1,1,0,0)
\end{eqnarray*}%
which corresponds to $[1,1]$. Similarly we have $a_{12}=a_{21}=(0,0,0,0)$
and $a_{22}=(1,1,0,0).$ Thus $B_2B_2^{-1}=Id.$

\section{Diagonalization}

\bigskip

\subsection{\protect\bigskip Eigenvalues and central eigenvalues}

Let $A$ be in $gl(n,\overline{\mathbb{IR}}).$ An eigenvalue of $A$ is an
element $\mathcal{X}\in \overline{\mathbb{IR}}$ such that there exists a
vector $\mathcal{V}\neq 0\in \overline{\mathbb{IR}}^{n}$ with
\begin{equation*}
A\cdot ^{t}\mathcal{V}=\mathcal{X}\cdot ^{t}\mathcal{V}.
\end{equation*}%
Thus $\mathcal{X}$ is a root of the characteristical polynomial with
coefficients in the ring $\overline{\mathbb{IR}}$%
\begin{equation*}
C_{A}(\mathcal{X})=\det (A-\mathcal{X}I)=0.
\end{equation*}%
\textbf{Example.} Let
\begin{equation*}
B_3=\left(
\begin{array}{cc}
\lbrack 1,2] & [1,2] \\
\lbrack 1,3] & [2,5]%
\end{array}%
\right) .
\end{equation*}%
We have
\begin{equation*}
B_3-\mathcal{X}I=\left(
\begin{array}{cc}
\lbrack 1,2]\smallsetminus \mathcal{X} & [1,2] \\
\lbrack 1,3] & [2,5]\smallsetminus \mathcal{X}%
\end{array}%
\right)
\end{equation*}%
and
\begin{eqnarray*}
\det (B_3-\mathcal{X}I) &=&([1,2]\smallsetminus \mathcal{X}%
)([2,5]\smallsetminus \mathcal{X})-[1,3][1,2] \\
&=&[2,10]-\mathcal{X}[2,5]-\mathcal{X}[1,2]+(\smallsetminus \mathcal{X}%
)(\smallsetminus \mathcal{X})-[1,6] \\
&=&(\smallsetminus \mathcal{X})(\smallsetminus \mathcal{X})-\mathcal{X}%
[3,7]+[1,4].
\end{eqnarray*}%
Let $\mathcal{X}=([x,y],0).$ It is represented in $\mathcal{A}_{4}$ by $%
(x,y,0,0)$ or $(0,y,x,0)$ or $(0,0,x,y)=-(x,y,0,0).$

\bigskip

\textbf{First case}: \ $\det (B_3-\mathcal{X}%
I)=(x^{2},y^{2},0,0)-(3x,7y,0,0)+(1,4,0,0)=(x^{2}-3x+1,y^{2}-7y+4,0,0).$
Then $\det (B_3-\mathcal{X}I)=0$ implies%
\begin{equation*}
\left\{
\begin{array}{c}
x^{2}-3x+1=0, \\
y^{2}-7y+4=0,%
\end{array}%
\right.
\end{equation*}%
that is
\begin{equation*}
\left\{
\begin{array}{c}
x=\dfrac{3\pm \sqrt{5}}{2}, \\
y=\dfrac{7\pm \sqrt{33}}{2}.%
\end{array}%
\right.
\end{equation*}%
We obtain
\begin{equation*}
\left\{
\begin{array}{c}
\mathcal{X}_{1}=([\dfrac{3+\sqrt{5}}{2},\dfrac{7+\sqrt{33}}{2}],0), \\
\mathcal{X}_{2}=([\dfrac{3-\sqrt{5}}{2},\dfrac{7+\sqrt{33}}{2}],0), \\
\mathcal{X}_{3}=([\dfrac{3-\sqrt{5}}{2},\dfrac{7-\sqrt{33}}{2}],0). \\
\end{array}%
\right.
\end{equation*}

\textbf{Second case: } $\det (B_3-\mathcal{X}%
I)=(0,y^{2}+x^{2},2xy,0)-(0,7y,7x,0)+(1,4,0,0)=(1,y^{2}+x^{2}-7y+4,2xy-7x,0).
$ Then $\det (B_3-\mathcal{X}I)=0$ implies%
\begin{equation*}
\left\{
\begin{array}{c}
1-2xy+7x=0, \\
y^{2}+x^{2}-7y+4=0.%
\end{array}%
\right.
\end{equation*}%
This gives
\begin{equation*}
4y^{4}-56y^{3}+261y^{2}-455y+197=0.
\end{equation*}%
We have the following solutions
\begin{equation*}
(x;y)=\{(-2,8;3,32),(2,9;3,67),(-0,17;0,63),(0,17;6,37)\}.
\end{equation*}%
We obtain the eigenvalues
\begin{equation*}
\left\{
\begin{array}{c}
\mathcal{X}_{4}=([-2,8,\ 3.32],0), \\
\mathcal{X}_{5}=([-0.17,\ 0.63],0). \\
\end{array}%
\right.
\end{equation*}

\qquad

\textbf{Third} \textbf{\ case}: $\det (B_3-\mathcal{X}%
I)=(x^{2},y^{2},0,0)+(3x,7y,0,0)+(1,4,0,0)=(x^{2}+3x+1,y^{2}+7y+4,0,0).$
Then $\det (B_3-\mathcal{X}I)=0$ implies%
\begin{equation*}
\left\{
\begin{array}{c}
x=\dfrac{-3\pm \sqrt{5}}{2}, \\
y=\dfrac{-7\pm \sqrt{33}}{2},%
\end{array}%
\right.
\end{equation*}%
then
\begin{equation*}
\begin{array}{c}
\mathcal{X}_{6}=([\dfrac{-3-\sqrt{5}}{2},\dfrac{-7+\sqrt{33}}{2}],0). \\
\end{array}%
\end{equation*}%
We obtain six eigenvalues.

\bigskip

\noindent\textbf{Remark. }To compute the interval-eigenvalues of a matrix $A$%
, we have to find the roots of the characteristical polynomial of $A$. But
this polynomial is with coefficients in $\overline{\mathbb{IR}}$ (or $%
\mathcal{A}_{4}$) and this set is not a field neither a factorial ring. Then
it is natural to meet some special results (e.g if we consider the second
degree polynomial $X^{2}-1$ with coefficients in $\frac{\mathbb{Z}}{8\mathbb{%
Z}}$ which is not factorial, it admits four roots,$1,3,5,7.$) In our example
we finds $6$ roots. Now if we consider the real matrix whose coefficients
are the centers of interval-coefficients of $B_3$, that is
\begin{equation*}
c_{B_3}=\left(
\begin{array}{cc}
1.5 & 1.5 \\
2 & 3.5%
\end{array}%
\right)
\end{equation*}%
then the eigenvalues of $c_{B_3}$ are $4.5$ and $0.5$ which are closed to
the center of $\mathcal{X}_{1}$ and $\mathcal{X}_{3}.$ We call these
eigenvalues, the central eigenvalues.

\begin{definition}
Let $A$ be a matrix in $gl(n,\overline{\mathbb{IR}}).$ Lat $A_{c}$ be the
real matrix whose elements are the center of the intervals of $A$. We say
that an eigenvalue of $A$ is a central eigenvalue if its center is (close
to) an eigenvalue of $A_{c}$.
\end{definition}

\bigskip \noindent\textbf{Remark.} The determination of negative eigenvalues
that is of type $(0,K)$ is similar. Nevertheless we have to consider only
matrices with positive entries thus we studies only the positive
eigenvalues. The negative eigenvalues do not correspond to physical entities.

\subsection{Eigenvectors, eigenspaces}

Now we will look the problem of reduction of an interval matrix. Recall that
the characteristical polynomial is with coefficient in a non factorial ring.
This is the bigest change with respect the classical real linear algebra.

\begin{defi}
Let $A$ a square matrix with coeffcients in $\overline{\mathbb{IR}}$. If $%
\mathcal{X}$ is an eigenvalue of $A$, then every vector $\mathcal{V}\in
\overline{\mathbb{IR}}^{n}$ satisfying $A^{t}\mathcal{V=X}^{t}\mathcal{V}$
is an eigenvector associated with $\mathcal{X}.$
\end{defi}

Let $E_{\mathcal{X}}$ be the set \
\begin{equation*}
E_{\mathcal{X}}=\{\mathcal{V}\in \overline{\mathbb{IR}}^{n}\text{ such that }%
A^{t}\mathcal{V=X}^{t}\mathcal{V}\}.
\end{equation*}%
Then $E_{\mathcal{X}}$ is a $\mathbb{R}$-subspace of $\overline{\mathbb{IR}}%
^{n}$ where $n$ is the order of the matrix $A$. It is also a $\overline{%
\mathbb{IR}}$ submodule of $\overline{\mathbb{IR}}^{n}$.

\begin{propo}
Let $\mathcal{X}_{1}$ and $\mathcal{X}_{2}$ be two distinguish eigenvalues
of $A $. Then $E_{\mathcal{X}_{1}}\cap E_{\mathcal{X}_{2}}=\{0\}.$
\end{propo}

\textit{Proof}\textbf{. }Let $\mathcal{V}$ be in $E_{\mathcal{X}_{1}}\cap E_{%
\mathcal{X}_{2}}$. We have
\begin{eqnarray*}
A^{t}\mathcal{V} &\mathcal{=}&\mathcal{X}_{1}\mathcal{V}, \\
A^{t}\mathcal{V} &\mathcal{=}&\mathcal{X}_{2}\mathcal{V}.
\end{eqnarray*}%
This $\mathcal{X}_{1}\mathcal{V}\smallsetminus \mathcal{X}_{2}\mathcal{V}=(%
\mathcal{X}_{1}\smallsetminus \mathcal{X}_{2})\mathcal{V}=0$. As $\overline{%
\mathbb{IR}}$ is without zero divisor, we have $\mathcal{X}%
_{1}\smallsetminus \mathcal{X}_{2}=0$ or $\mathcal{V}=0$. We deduce $E_{%
\mathcal{X}_{1}}\cap E_{\mathcal{X}_{2}}=\{0\}.$

\bigskip

\begin{propo}
Let $C_{A}(\mathcal{X})$ be the characteristical polynomial of $A$. If the
real polynomial $C_{C_{A}}(\mathcal{X})$ associated with the central matrix
of $A$ is a product of factor of degree $1$, then $C_{A}(\mathcal{X})$
admits a factorization on $\overline{\mathbb{IR}}$
\end{propo}

We have seen that $C_{A}(\mathcal{X})$ can be have more than degree($C_{A}(%
\mathcal{X}$) roots. If $\mathcal{X}_{1},\cdots ,\mathcal{X}_{n}$ are the
central roots, we have the decomposition
\begin{equation*}
C_{A}(\mathcal{X})=a_{n}\prod\limits_{i=1}^{n}(\mathcal{X}\smallsetminus
\mathcal{X}_{i}).
\end{equation*}

\medskip

\noindent\textbf{Example. }If we consider the matrix
\begin{equation*}
B_3=\left(
\begin{array}{cc}
\lbrack 1,2] & [1,2] \\
\lbrack 1,3] & [2,5]%
\end{array}%
\right) .
\end{equation*}%
then $C_{B_3}(\mathcal{X})$ admits $\mathcal{X}_{1},\cdots ,\mathcal{X}_{6}$
as positive roots. The central eigenvalues are $\mathcal{X}_{1}$ and $%
\mathcal{X}_{3}$ and we have%
\begin{equation*}
\det (B_3\smallsetminus \mathcal{X}I)=(\mathcal{X}\smallsetminus \mathcal{X}%
_{1})(\mathcal{X}\smallsetminus \mathcal{X}_{3}).
\end{equation*}%
If we consider the roots $\mathcal{X}_{2}=([\dfrac{3-\sqrt{5}}{2},\dfrac{7+%
\sqrt{33}}{2}],0)$, and if we assume that $C_{B_3s}(\mathcal{X})=(\mathcal{X}%
\smallsetminus \mathcal{X}_{2})(\mathcal{X}\smallsetminus Y)$, we obtain
\begin{equation*}
Y=(3,7,\frac{3-\sqrt{5}}{2},\frac{7+\sqrt{33}}{2})=(\frac{3+\sqrt{5}}{2},%
\frac{7-\sqrt{33}}{2},0,0)
\end{equation*}%
which does not correspond to a positive eigenvalue.

\bigskip

\begin{theorem}
For any $n$-ule of roots $(\mathcal{X}_{1},\cdots ,\mathcal{X}_{n})$ such
that $C_{A}(\mathcal{X})=a_{n}\prod\limits_{i=1}^{n}(\mathcal{X}%
\smallsetminus \mathcal{X}_{i})$, and if for any $i=1,\cdots ,n$ the
dimension of $E_{\mathcal{X}_{i}}$ coincides with the multiplicity of $%
\mathcal{X}_{i},$ then we have the vectorial decomposition $\overline{%
\mathbb{IR}}^{n}=\oplus _{i\in I}E_{\mathcal{X}_{i}}$ where the roots $%
\mathcal{X}_{i},i\in I$ are pairwise distinguish.
\end{theorem}

\bigskip

\noindent\textbf{Example.} Let us compute the eigenspaces of $B_3$
associated to the central eigenvalues.

\begin{itemize}
\item $\mathcal{X}_{1}=([\dfrac{3+\sqrt{5}}{2},\dfrac{7+\sqrt{33}}{2}],0).$

\noindent Let $V=\left(
\begin{array}{c}
V_{1} \\
V_{2}%
\end{array}%
\right) \in \overline{\mathbb{IR}}.$ Then
\begin{equation*}
(A-\mathcal{X}_{1}I)\left(
\begin{array}{c}
V_{1} \\
V_{2}%
\end{array}%
\right) =0
\end{equation*}%
is equivalent to
\begin{equation*}
\left(
\begin{array}{ll}
(0,[\dfrac{1+\sqrt{5}}{2},\dfrac{3+\sqrt{33}}{2}])V_{1} & ([1,2],0)V_{2} \\
([1.3],0)V_{1} & (0,[\dfrac{-1+\sqrt{5}}{2},\dfrac{-3+\sqrt{33}}{2}])V_{2}%
\end{array}%
\right) =0
\end{equation*}%
that is
\begin{equation*}
\left\{
\begin{array}{c}
\smallsetminus \lbrack \dfrac{1+\sqrt{5}}{2},\dfrac{3+\sqrt{33}}{2}%
]V_{1}+[1,2]V_{2}=0, \\
\lbrack 1.3]V_{1}\smallsetminus \lbrack \dfrac{-1+\sqrt{5}}{2},\dfrac{-3+%
\sqrt{33}}{2}]V_{2}=0.%
\end{array}%
\right.
\end{equation*}%
This gives
\begin{equation*}
V_{2}=\frac{[\dfrac{1+\sqrt{5}}{2},\dfrac{3+\sqrt{33}}{2}]V_{1}}{[1,2]}.
\end{equation*}%
If we choose $V_{1}=([1,1],0)$ we have
\begin{equation*}
\begin{array}{llll}
\medskip V_{2} & = & \frac{[\dfrac{1+\sqrt{5}}{2},\dfrac{3+\sqrt{33}}{2}]}{%
[1,2]} &  \\
\medskip  & = & [\dfrac{-1+\sqrt{5}}{2},\dfrac{-3+\sqrt{33}}{2}]\bullet
(\smallsetminus \lbrack -1,\frac{-1}{2}]) &  \\
& = & \smallsetminus ([-\frac{3+\sqrt{33}}{2},-\frac{1+\sqrt{5}}{4}]) &
\end{array}%
\end{equation*}%
Thus the $\mathcal{X}_{1}$-eigenvectors are of the form
\begin{equation*}
V=\left(
\begin{array}{c}
([1,1],0) \\
\smallsetminus ([-\frac{3+\sqrt{33}}{2},-\frac{1+\sqrt{5}}{4}],0)%
\end{array}%
\right) .
\end{equation*}
\end{itemize}

\noindent\textbf{Remark. }We can choose $V_{1}$ such that all tje
coordinatez of $V$ are positive. For example if $V_{1}=[1,2]$ then $V=\left(
\begin{array}{c}
([1,2],0) \\
([\dfrac{1+\sqrt{5}}{2},\dfrac{3+\sqrt{33}}{2}],0)%
\end{array}%
\right) $

\begin{itemize}
\item $\mathcal{X}_{3}=([\dfrac{3-\sqrt{5}}{2},\dfrac{7-\sqrt{33}}{2}],0).$

\noindent Let $V=\left(
\begin{array}{c}
V_{1} \\
V_{2}%
\end{array}%
\right) \in \overline{\mathbb{IR}}.$ Then
\begin{equation*}
(A-\mathcal{X}_{1}I)\left(
\begin{array}{c}
V_{1} \\
V_{2}%
\end{array}%
\right) =0
\end{equation*}%
is equivalent to
\begin{equation*}
\left(
\begin{array}{ll}
([\dfrac{-1+\sqrt{5}}{2},\dfrac{-3+\sqrt{33}}{2}],0)V_{1} & ([1,2],0)V_{2}
\\
([1.3],0)V_{1} & ([\dfrac{1+\sqrt{5}}{2},\dfrac{3+\sqrt{33}}{2}],0)V_{2}%
\end{array}%
\right) =0
\end{equation*}%
that is
\begin{equation*}
\left\{
\begin{array}{c}
\lbrack \dfrac{-1+\sqrt{5}}{2},\dfrac{-3+\sqrt{33}}{2}]V_{1}+[1,2]V_{2}=0,
\\
\lbrack 1.3]V_{1}+[\dfrac{1+\sqrt{5}}{2},\dfrac{-3+\sqrt{33}}{2}]V_{2}=0.%
\end{array}%
\right.
\end{equation*}%
This gives
\begin{equation*}
V_{2}=\frac{\smallsetminus ([\dfrac{-1+\sqrt{5}}{2},\dfrac{-3+\sqrt{33}}{2}%
])V_{1}}{[1,2]}.
\end{equation*}%
If we choose $V_{1}=([1,1],0)$ we have
\begin{equation*}
\begin{array}{llll}
\medskip V_{2} & = & \frac{\smallsetminus \lbrack \dfrac{-1+\sqrt{5}}{2},%
\dfrac{-3+\sqrt{33}}{2}]}{[1,2]} &  \\
\medskip  & = & \smallsetminus \lbrack \dfrac{-1+\sqrt{5}}{2},\dfrac{-3+%
\sqrt{33}}{2}]\bullet (\smallsetminus \lbrack -1,\frac{-1}{2}]) &  \\
& = & ([\frac{3-\sqrt{33}}{2},\frac{1-\sqrt{5}}{4}]) &
\end{array}%
\end{equation*}%
Thus the $\mathcal{X}_{3}$-eigenvectors are of the form
\begin{equation*}
V=\left(
\begin{array}{c}
([1,1],0) \\
([\frac{3-\sqrt{33}}{2},\frac{1-\sqrt{5}}{4}],0)%
\end{array}%
\right) .
\end{equation*}
\end{itemize}

\section{The Exponential map}

We define the exponential map
$$Exp :gl(n,\overline{\mathbb{IR}})\longrightarrow gl(n,\overline{\mathbb{IR}}) $$
in a classical way by series expansions. If the matrix $A$ is diagonalizable, then
$$D=P^{-1}AP$$
is diagonal and $Exp(A)$ is a diagonal matrix whose diagonal element are the exponential of the eigenvalues.

\noindent {\bf Example.} Let
\begin{equation*}
B_3=\left(
\begin{array}{cc}
\lbrack 1,2] & [1,2] \\
\lbrack 1,3] & [2,5]%
\end{array}%
\right) .
\end{equation*}%
The central eigenvalues are
\begin{equation*}
\left\{
\begin{array}{c}
\mathcal{X}_{1} =([\dfrac{3+\sqrt{5}}{2},\dfrac{7+\sqrt{33}}{2}],0), \\
\mathcal{X}_{3}=([\dfrac{3-\sqrt{5}}{2},\dfrac{7-\sqrt{33}}{2}],0). \\
\end{array}%
\right.
\end{equation*} and we have
\begin{equation*}
D=\left(
\begin{array}{cc}
\mathcal{X}_{1} & 0 \\
0 & \mathcal{X}_{3}
\end{array}%
\right)
\end{equation*}%
with
\begin{equation*}
P=\left(
\begin{array}{cc}
([1,1],0) & ([1,1],0) \\
\smallsetminus ([-\frac{3+\sqrt{33}}{2},-\frac{1+\sqrt{5}}{4}],0)& ([\frac{3-\sqrt{33}}{2},\frac{1-\sqrt{5}}{4}],0)
\end{array}%
\right) .
\end{equation*}%
We deduce
$$Exp(B_3)=P.\left(
\begin{array}{cc}
 ([exp(\dfrac{3+\sqrt{5}}{2}),exp(\dfrac{7+\sqrt{33}}{2}],0)) & 0 \\
0 & ([exp(\dfrac{3-\sqrt{5}}{2}),exp(\dfrac{7-\sqrt{33}}{2}],0))
\end{array}%
\right).P^{-1}$$
In a forthcomming paper, we apply this calculus to solve linear differential system.


\begin{thebibliography}{9}
\bibitem{GN-R} Goze Nicolas, Remm Elisabeth. An algebraic approach to the
set of intervals. arXiv math 0809.5150 (2008)
\end{thebibliography}
\end{document}